\documentclass[11pt,a4paper,twoside]{paper}
\newif\ifRR\RRfalse
\ifRR
  \usepackage{RR}
\fi

\usepackage[T1]{fontenc}
\usepackage{amsmath,amssymb,amsthm}
\usepackage{parskip}
\newtheoremstyle{etoile}{\parskip}{\parskip}{\itshape}
                        {0pt}{\bfseries\sffamily}{.}{ }{}
\theoremstyle{etoile}

\newcommand\bbox{\hfill\rule{2mm}{2mm}}

\newtheorem{prop}{Proposition}[section]

\newtheorem{theo}[prop]{Theorem}
\newtheorem{defin}[prop]{Definition}
\newtheorem{lem}[prop]{Lemma}

\makeatletter
\@addtoreset{equation}{section}
\makeatother

\newcommand\egaldef{\stackrel{\mbox{\upshape\tiny def}}{=}}

\newcommand\1{\leavevmode\hbox{\rm \small1\kern-0.35em\normalsize1}}
        
\newcommand\EE{\mathsf{E}}
\newcommand\Dc{\mathcal{D}} 
\newcommand\Gb{\mathbf{G}}
\newcommand\Ac{\mathcal{A}}
\newcommand\Bc{\mathcal{B}}
\newcommand\Cc{\mathcal{C}}
\newcommand\Lc{\mathcal{L}}
\newcommand\Mc{\mathcal{M}}
\newcommand\Oc{\mathcal{O}}
\newcommand\Rc{\mathcal{R}}
\newcommand\Uc{\mathcal{U}}
\newcommand\Sb{\mathbf{S}}
\newcommand\ZZ{\mathbb{Z}}
\newcommand\n{^{\scriptscriptstyle (N)}}
\newcommand\nn{\scriptscriptstyle N}
\def\DD{\displaystyle} 
\DeclareMathOperator*{\lra}{\leftrightarrows}
\DeclareMathOperator*{\rla}{\rightleftarrows}

\begin{document}
 
\ifRR
\else
\title{Stochastic Dynamics of Discrete Curves and Exclusion Processes. \\[0.2cm]
Part 1:  Hydrodynamic Limit of the ASEP System}

\author{Guy Fayolle \thanks{INRIA - Domaine de Voluceau, 
Rocquencourt
BP 105 - 78153 Le Chesnay Cedex - France. Contact:
\texttt{Guy.Fayolle@inria.fr, Cyril.Furtlehner@inria.fr}} 
\and Cyril Furtlehner \footnotemark[1]} 
\date{December 2005}

\maketitle
\begin{abstract}
This preliminary report is the foreword of a series dedicated to stochastic deformations of curves. Problems  are set in terms of exclusion processes, the ultimate goal being to derive  hydrodynamic limits for these systems after proper scalings. Here, only the basic \textsc{asep} system on the torus  $\ZZ/N\ZZ$ is analyzed.  The usual sequence of empirical measures,  converges in probability to a deterministic measure, which is the unique weak solution of a Cauchy problem. The method  presents some new features,  letting hope for extensions to higher dimension. It relies on the analysis of specific partial differential equations involving variational calculus. Namely, the variables are the values of functions at given points and their number becomes infinite. 
\end{abstract}


\fi


\section{Preliminaries}\label{INTRO}
 Interplay between discrete and continuous  description is a recurrent   question
in statistical physics, which  in some cases can be answered quite rigorously via probabilistic methods. In the context of reaction-diffusion 
systems, this is tantamount to studying  fluid or hydrodynamics limits.   
Number of approaches have been proposed, in particular in the framework of exclusion 
processes, see e.g. \cite{Li},\cite{MaPr} \cite{Sp}, \cite{KiLa} 
and references therein.  As far as fluid limits are at stake, all these methods 
have in common to be limited to systems for which the stationary states are  
given in closed product forms, as far as hydrodynamic limits are concerned, 
or at least for which the invariant measure for finite $N$ (the size of the system) is explicitly known. For instance,  \textsc{asep} with open boundary can be described in terms of matrix product form (a sort of a non-commutative product form) and the continuous  limits can be understood by means of brownian bridges (see 
\cite{DeEvHaPa}). We propose to adress this question from the following different point of view: starting from discrete sample paths subjected to stochastic deformations, the ultimate goal is to understand the nature of the limit curves when $N$ increases to infinity. How do these curves evolve with time, and which limiting process do they represent as $t$ goes to infinity (equilibrium curves)? Following \cite{FaFu} and 
\cite{FaFu2},  we will try to give some  partial answers  to these questions in a series of papers.

This first  study is mainly dedicated to the \textsc{asep} model. The mathematical approch relies on the analysis of specific partial differential equations involving variational calculus. A usual sequence of empirical measures is shown to converge in probability to a deterministic measure, which is the unique weak solution of a Cauchy problem. Here variables are the values of some function at given points and their number becomes infinite.  In our opinion, the method  presents some new features, which let hope for extensions to higher dimension. 
  
A future concern will be to establish a complete  hierarchy of systems of hydrodynamic equations , whose steady state will help to describe non-Gibbs states.

\section{Model definition}
\subsection{A stochastic clock model}\label{sec:clock}
The systems we will consider  can typically be describded as an oriented path embedded in a bidimensional manifold:  $N$ steps of equal size, each one being chosen among a discrete set of $n$ possible orientations drawn from the set $\{2k\pi/n, k=0\ldots n-1\}$ of angles with some given origin. The stochastic dynamics which is applied consists in displacing one single point at a time without breaking the path, while keeping 
all links within the set of admissible orientations. In this operation
two links are displaced. This constrains quite strongly the possible dynamical rules, 
which are given in terms of \emph{reactions}  between consecutive links. 

We have
\[
X^k X^l\ \rla_{\lambda^{lk}}^{\lambda^{kl}} \ X^l X^k,\quad
k\in[1, n],\, k \ne l .
\]

\subsection{Examples} \label{sec:exemple}
\emph{(1) The simple exclusion process}

The first elementary and most sudied example is the simple exclusion process: this model, after mapping particles onto links, corresponds to a one-dimensional fluctuating interface. Here we have a binary alphabet and letting
$X^1=\tau$ and $X^2=\bar\tau$, the set of reactions simply rewrites
\[
\tau\bar\tau \lra_{\lambda^+}^{\lambda^-}\ \bar\tau\tau,
\]
where $\lambda^\pm$ are the transition rates for the jump of a particle to
the right or to the left. 

\emph{(2) The triangular lattice and the ABC model}

 Here the evolution of the random walk is  restricted to the triangular
 lattice. Each link (or step) of the walk is either $1$, $e^{2i\pi/3}$
 or $e^{4i\pi/3}$, and quite naturally will be said to be of type A, B or C, respectively. 
 This corresponds to the so-called \emph{ABC
 model}, since there is a coding by means of a  $3$-letter alphabet. 
 The set of \emph{transitions} (or reactions) is given by
\begin{eqnarray}
AB\ \lra_{p^+}^{p^-}\ BA, \qquad 
BC\ \lra_{q^+}^{q^-}\ CB, \qquad 
CA\ \lra_{r^+}^{r^-}\ AC, \qquad 
\end{eqnarray}
where there is a priori no symmetry, but we will impose \emph{periodic
boundary conditions} on the sample paths. This model was first
introduced in \cite{EvFoGoMu} in the context of particles with exclusion, and for some cases corresponding to the reversibility of the process, a Gibbs form for the invariant measure was given in \cite{EvKaKoMu}

\section{Hydrodynamics for a one-dimensional asymmetric exclusion process 
[ASEP]}

As mentioned above, we aim at obtaining hydrodynamic equations for a  class of exclusion  models. The method, although relying on classical powerful tools  (martingales, relative compactness of measures, functional analysis), has some new features which should hopefully prove fruitful in other contexts. The essence of the approach is in fact contained in the analysis of the popular \textsc{asep} model, presented below. We note the difficulty to find in the existing literature a complete study encompassing  various special cases (symmetry, total asymmetry, etc). Some proofs will only be sketched, and the related results  presented as claims or even conjectures.

Consider $N$ sites labelled from $1$ to $N$, forming a discrete closed curve in the plane, so that  the numbering of sites is implicitly  taken modulo $N$, i.e. on the discrete torus $\Gb\n\egaldef \ZZ/N\ZZ$.  In higher dimension, say on the lattice $\ZZ^k$, the related  set of sites would be drawn on the torus  $(\ZZ/N\ZZ)^k$. 

We gather below some notational material valid throughout this section.
\begin{itemize}
\item $\Rc$ stands for the real line. $\Cc^k[0,1]$ is the collection of all real-valued, $k$-continuously differentiable functions defined on the interval $[0,1]$, and $\Mc$  is the space of all  finite positive measures on the torus $\Gb\egaldef [0,1)$. 

$\Cc^\infty_{0}(K)$  is the space of infinitely differentiable functions with compact support
 included in $K$.
\item For $\Sb$ an arbitrary metric space, $\mathcal{P}(\Sb)$ is the set of probability measures $\Dc_\Sb[0,T]$ is the space of right continuous functions $z : [0,\infty]\to \Sb$ with left limits and $t\to z_{t}$.
\item For $ i=1,\ldots,N$, let  $A_i\n(t)$ and $B_i\n(t)$ be binary random variables representing respectively a particule or a hole at site $i$, so that, owing to the exclusion constraint,  $A_i\n(t)+B_i\n(t)=1$, for all $1\le i\le N$. Thus $\bigl\{\mathbf{A}\n(t)\egaldef\bigl(A_i\n(t), \ldots,A_N\n(t)\bigr), t\ge0\bigr\}$ is a Markov  process.
\item  $\Omega\n$ will denote the generator of the Markov process $\mathbf{A}\n(t)$, and $\mathcal{F}_t\n=\sigma\bigl(\mathbf{A}\n(s),s\leq t\bigr)$ is the associated  natural filtration.
\item Our purpose is to  analyze the sequence of empirical random measures 
\begin{equation} \label{eq:emp1}
\mu\n_{t} = \frac{1}{N}\sum_{i\in\Gb\n} A_i\n(t) \delta _\frac{i}{\nn},
\end{equation}
when $N\to\infty$, after a convenient scaling of the parameters of the generator 
$\Omega\n$. The probability distribution associated with the path of the Markov process
$\mu\n_t, t\in[0,T]$, for some fixed $T$, is simply denoted by $Q\n$.
 \end{itemize}

As usual, one can embed $\Gb\n$ in  $\Gb$, so that  a point  $i\in \Gb\n$ corresponds to the point  $i/N$ in $\Gb$. Hence, in view of (\ref{eq:emp1}), it is quite natural to let the sequence $Q\n$ be defined on a unique space  $\Dc_\Mc[0,T]$, which becomes a polish space (i.e. complete and separable) via the usual Skorokod topology, as soon as $\Mc$ is itself Polish (see e.g. \cite{EtKu}, chapter $4$). Without further comment,  $\Mc$ is assumed to be endowed with the vague product topology, as a consequence of the famous  Banach-Alaoglo and Tychonoff  theorems (see e.g. \cite{RUD,Ka}). 

Let  $\phi_a, \phi_b$ be two arbitrary functions in $\mathbf{C}^2[0,1]$ and
define the real-valued positive measure
\begin{equation} \label{eq:emp2}
Z\n _t[\phi_a,\phi_b] 
\egaldef \exp\biggl[\frac{1}{N}\sum_{i\in\Gb\n} \phi_a\bigl(\frac{i}{N}\bigr)A_i\n(t)+\phi_b\bigl(\frac{i}{N}\bigr)B_i\n(t) \biggr],
\end{equation} 
which is a  functional of $\phi_a,\phi_b$. For the sake of brevity, the explicit dependence of $A_i\n(t), B_i\n(t), Z\n_t[\phi_a,\phi_b]$  on $N,t,\phi$, will be omitted wherever the meaning remains clear from the context: for instance, we often shall simply write $A_i, B_i$ or $Z\n_t$. Also $Z\n$ stands for the process $\{Z\n_t, \,t\ge 0\}$.

A standard powerful method to prove the convergence (in a sense to be specified later) of  the sequence of probability measures introduced in (\ref{eq:emp1}) consists  first in showing its relative compactness, and then in verifying  the coincidence of all possible limit points (see e.g. \cite{KiLa}. Moreover here it suffices to prove these two properties for the sequence of projected measures defined on  $\Dc_\Rc[0,T]$ and corresponding to the processes $\{Z\n_t[\phi_a,\phi_b], t\ge 0\}$,  since the functions $\phi$ belong to 
 $\mathbf{C}^2[0,1]$.

Let us now introduce  quantities which, as far as scaling is concerned, are crucial
in order to obtain meaningful hydrodynamic equations.
\begin{equation}\label{eq:scale1}
\begin{cases}
\DD\lambda(N) \egaldef \frac{\lambda_{ab}(N)+\lambda_{ba}(N)}{2} ,\\[0.2cm]
\mu(N)  \egaldef \lambda_{ab}(N)-\lambda_{ba}(N),
\end{cases}
\end{equation}
where the dependence of the rates on $N$ is explicitly mentioned.

\begin{theo} \label{theo:main} Let system (\ref{eq:scale1})  have a given asymptotic expansion of the form
\begin{equation}\label{eq:scale2}
\begin{cases}
\DD\lambda(N) \egaldef \lambda N^2 +o(N^2) ,\\[0.2cm]
\mu(N)  \egaldef \mu N + o(N),
\end{cases}
\end{equation}
 where  $\lambda$ and $\mu$ are fixed positive constants. [As for the scaling assumption (\ref{eq:scale2}),  the random measure $\log Z\n_t$ is a functional of the underlying Markov process, in which the time has been speeded up by a factor $N^2$ and the space shrunk by $N^{-1}$]. Assume moreover the sequence of initial empirical measures $\log Z\n _0$, taken at time $t=0$, converges in probability to some deterministic measure with a given density $\rho(x,0)$, so that
  \begin{equation}\label{eq:init}
  \lim_{N\to\infty} \log Z\n _0 = \int_0^1[\rho(x,0)\phi_a(x)+(1-\rho(x,0))\phi_b(x)]dx, \quad\textrm{in probability},
 \end{equation}
for any pair of  functions $\phi_a, \phi_b \in\Cc^\infty_{0}(K)$, where $K\in\Rc$ is a compact containing the interval $[0,1]$. 

Then,  for every $t>0$, the sequence of random measures $\mu\n_t$ converges in probabilitys, as  $N\to\infty$, to a deterministic measure having a density $(\rho(x,t)$ with respect to the Lebesgue measure, which is the unique \emph{weak solution of  the   Cauchy problem}

\begin{equation}\label{eq:Cauchy}
\begin{split}
 \int_0^T \int_0^1 \left[\rho(x,t) \Bigl(\frac{\partial\theta(x,t)}{dt} +\lambda\frac{\partial^2\theta(x,t)}{dx^2}\Bigr) -\mu\rho(x,t)\bigl(1-\rho(x,t)\bigr)\frac{\partial\theta(x,t)}{dx}\right]dx dt \\[0.2cm]
 =  \int_0^1 \left[\rho(x,T)\theta(x,T) - \rho(x,0)\theta(x,0)\right] dx ,
\end{split}
 \end{equation}
where  (\ref{eq:Cauchy}) holds for any function $\theta\in\Cc^\infty_{0}([0,1]\times[0,T])$. 

 If, moreover, one assumes the existence of $\frac{\partial^2\rho(x,0)}{dx^2}$, then (\ref{eq:Cauchy}) reduces to a classical Burger's equation
\[
\frac{\partial\rho(x,t)}{dt} = \lambda\frac{\partial^2\rho(x,t)}{dx^2} +
\mu [1-2\rho(x,t)] \frac{\partial\rho(x,t)}{dx}.
\]
\end{theo}

\begin{proof} \mbox{ } The sketch of the proof  is spread over three main subsections, referred to hereafter  as \textbf{P1}, \textbf{P2} and \textbf{P3}.

\paragraph{P1 [Existence of limit points: sequential compactness]}
As usual in problems dealing with convergence of sequences of probability measures, our very starting point will be to establish the weak relative compactness of the set $\{\log Z\n_t, N\ge1\}$. Some of the probabilistic arguments employed in this paragraph are in a way  classical and can be found in good books, e.g. \cite{Sp, KiLa}, although for simpler models.

The process 
\begin{equation}\label{eq:martin1}
U_t\n \egaldef  Z\n_t - Z\n_0 - \int_0^t \Omega\n  [Z\n_s] ds 
\end{equation}
is a bounded $\{\mathcal{F}_t\n\}$-martingale. Using the exponential form of $Z\n_t$ together with classical stochastic calculus (see e.g. \cite{EtKu}, chap.3, page 93), it follows that 
\begin{equation}\label{eq:martin2}
[V_t\n]  \egaldef (U\n_t)^2  - \int_0^t \Bigl(\Omega\n  [(Z\n_s)^2]
- 2 Z\n_s \Omega\n  [Z\n_s] \Bigr)ds 
\end{equation}
is also a bounded real martingale. 

From $A_i\n(t)+B_i\n(t)=1,\forall 1\leq i\leq N$, on sees that 
$Z\n _t$ is mainly a functional of the sole function  $\psi_{xy}\egaldef\phi_x-\phi_y=-\psi_{yx}$, up to a constant uniformly bounded in $N$. Hence, setting
\begin{eqnarray*}\label{eq:delta}
\Delta\psi_{xy} \Bigl(\frac{i}{N}\Bigr) &\egaldef & \psi_{xy}\Bigl(\frac{i+1}{N}\Bigr) -\psi_{xy}\Bigl(\frac{i}{N}\Bigr),  \\
 \widetilde{\lambda}_{xy}(i,N)& \egaldef & \lambda_{xy}(N)\left[
 \exp\biggl(\frac{1}{N} \Delta\psi_{xy} \Bigl(\frac{i}{N}\Bigr)\biggr) - 1\right], 
\quad \ xy = ab \ \text{or} \ ba , 
\end{eqnarray*}
we have
\begin{equation}\label{eq:gen1}
 \Omega\n [Z\n_t ]  =    L\n_t Z\n_t ,
 \end{equation}
 where
\begin{equation}\label{eq:gen2}
L\n_t  =  \sum_{i\in\Gb\n} \widetilde{\lambda}_{ab}(i,N) A_{i}B_{i+1} +  
 \widetilde{\lambda}_{ba}(i,N) B_{i}A_{i+1}.
\end{equation}
  
 By using the exclusion property,  a straightforward calculation in equation (\ref{eq:gen2}) allows to rewrite  (\ref{eq:martin2}) in  the form
\begin{equation}\label{eq:martin3}
[V_t\n] = (U\n_t)^2  - \int_0^t (Z\n_s)^2R_s\n ds ,
\end{equation}
where the process $R_t\n$ is stricly positive and given by
\[
R_t\n =  \sum_{i\in\Gb\n} \frac{[\widetilde{\lambda}_{ab}(i,N)]^2}{\lambda_{ab}(N)}
A_{i}B_{i+1} \,+\, \frac{[\widetilde{\lambda}_{ba}(i,N)]^2}{\lambda_{ba}(N)} B_{i}A_{i+1}.
\]  
The integral term in  (\ref{eq:martin3}) is nothing else but the increasing process associated with Doob's decomposition of the submartingale  $(U\n_t)^2$. 

The folllowing estimates are crucial.
\begin{lem} \label{lem:esti}
\begin{eqnarray}
L\n_t  & = & \Oc(1) ,  \label{eq:lem1}  \\[0.2cm]
R_t\n & = & \Oc\Bigl(\frac{1}{N}\Bigr). \label{eq:lem2}
\end{eqnarray}
\end{lem}
\begin{proof} We will derive (\ref{eq:lem1}) by estimating  the right-hand side member of equation (\ref{eq:gen2}). 

Clearly,  $\Delta\psi_{xy} \Bigl(\frac{i}{N}\Bigr) = \frac{1}{N}\psi_{xy} '\Bigl(\frac{i}{N}\Bigr) + \mathcal{O}\Bigl( \frac{1}{N^2}\Bigr)$,  where $\psi'$ denotes the derivative of $\psi$.  Then, taking a  second order expansion of the exponential function and using definitions (\ref{eq:scale1}) and (\ref{eq:scale2}), we can rewrite  (\ref{eq:gen2}) as
\begin{align}\label{eq:estim1}
 L\n_t & =  \frac{\mu(N)}{N} \sum_{i\in\Gb\n}\left[  \frac{A_{i}+A_{i+1}}{2} - A_{i}A_{i+1} \right] \Delta\psi_{ab} \Bigl(\frac{i}{N}\Bigr) \nonumber \\
&+  \frac{\lambda(N)}{N}\sum_{i\in\Gb\n} (A_{i}-A_{i+1})\Delta\psi_{ab} \Bigl(\frac{i}{N}\Bigr)
 +  \Oc \Bigl(\frac{1}{N}\Bigr).
\end{align}

 The first  sum in (\ref{eq:estim1}) is uniformly bounded by a  constant depending on 
 $\psi$. Indeed,  $|A_{i}| \le 1$ and $\psi\in\mathbf{C}^2[0,1]$,  so that  $\psi'$ is of bounded variation.
 
  As for the second sum coming in  (\ref{eq:estim1}), we have
  \[ 
  \sum_{i\in\Gb\n} (A_{i}-A_{i+1})\Delta\psi_{ab} \Bigl(\frac{i}{N}\Bigr) = 
   \sum_{i\in\Gb\n} A_{i+1} \left[\Delta\psi_{ab} \Bigl(\frac{i+1}{N}\Bigr) -
   \Delta\psi_{ab} \Bigl(\frac{i}{N}\Bigr)\right].
   \]
Then the discrete Laplacian 
\[\Delta\psi_{ab} \Bigl(\frac{i+1}{N}\Bigr) - \Delta\psi_{ab} \Bigl(\frac{i}{N}\Bigr) \equiv 
\psi_{ab}\Bigl(\frac{i+2}{N}\Bigr) -2\psi_{ab} \Bigl(\frac{i+1}{N}\Bigr) + \psi_{ab} \Bigl(\frac{i}{N}\Bigr)
\]
admits of the simple form
\begin{equation}\label{eq:psi}
\Delta\psi_{ab} \Bigl(\frac{i+1}{N}\Bigr) - \Delta\psi_{ab} \Bigl(\frac{i}{N}\Bigr) =
\frac{1}{N^2} \psi_{ab}''\Bigl(\frac{i}{N}\Bigr) +\Oc\Bigl(\frac{1}{N^2}\Bigr), 
\end{equation}
where $ \psi''$ denotes the second derivative of $\psi$. 

By  (\ref{eq:scale2}),  $\DD\lambda(N)=\lambda N^2 +o(N^2)$,  
so that (\ref{eq:psi}) implies
 \begin{equation}\label{eq:estim2}
 \frac{\lambda(N)}{N}\sum_{i\in\Gb\n} (A_{i}-A_{i+1})\Delta\psi_{ab} \Bigl(\frac{i}{N}\Bigr)
= \sum_{i\in\Gb\n} \frac{\lambda A_{i+1}}{N} \psi_{ab}''\Bigl(\frac{i}{N}\Bigr) + 
o\Bigl(\frac{1}{N}\Bigr) = \Oc(1),
\end{equation}
which concludes the proof of (\ref{eq:lem1}). The computation of $R_t\n$ leading to (\ref{eq:lem2}) can be obtained via similar arguments.
\end{proof}

To show the relative compactness of the family $Z\n$, which here, by separability and completeness of the underlying spaces, is equivalent to tightness, we proceed as in \cite{KiLa} by means of the following useful criterion.

\begin{prop} [Aldous's tightness criterion, see \cite{Bil}] \label{prop:tight}

A sequence $\{X\n\}$ of random elements of  $\Dc_\Rc[0,T]$ is tight (i.e. the distributions of the  $\{X\n\}$ are tight) if the two following conditions hold:

\begin{itemize} 
\item[(i) ]
\begin{equation} \label{C1}
\lim_{a\to\infty} \limsup_N P [ || X\n ||  \ge a] =0, 
\end{equation}
where $\DD || X\n || \egaldef \sup_{t\le T} |X\n_t|$.
\item[(ii)] For each $\epsilon, \eta$, there exist a $\delta_0$ and an $N_0$, such that, if 
$\delta\le\delta_0$ and $N\ge N_0$, and if $\tau$ is an arbitrary stopping time with
$\tau+\delta\le T$, then
\begin{equation} \label{C2}
    P\bigl[ | X\n_{\tau+\delta} - X\n_{\tau} | \ge \epsilon \bigr] \le \eta .
\end{equation}
\end{itemize}
Note that condition (\ref{C1}) is always necessary for tightness. \bbox
\end{prop}

We will now apply Lemma \ref{lem:esti} to equations (\ref{eq:martin1}) and  (\ref{eq:martin3}), the role of $X\n$ in Proposition \ref{prop:tight} being played by 
$Z\n$. 

Observe that, by the uniform boundedness  of $Z\n_t$, 
condition (\ref{C1}) is immediately verified.

To check condition (\ref{C2}), rewrite (\ref{eq:martin1}) as
\begin{equation} \label{eq:martin4}
Z\n_{t+\delta} - Z\n_t = U\n_{t+\delta} - U\n_t + \int_{t}^{t+\delta}\Omega\n  [Z\n_s] ds .
\end{equation}

The integral term in (\ref{eq:martin4}) is bounded in modulus by $K\delta$ [where $K$ is a constant uniformly bounded in $N$ and $\psi$] and hence satisfies (\ref{C2}).
We are left with the analysis of $U\n_t$. But, from  (\ref{eq:martin3}), (\ref{eq:lem2}) and Doob's inequality for sub-martingales, we have
\begin{eqnarray}
\EE\bigl[(U\n_{t+\delta} - U\n_t )^2\bigr] &=& \EE\left[\int_t^{t+\delta} (Z\n_s)^2R_s\n ds\right]  \le  \frac{C}{N}  ,  \nonumber \\[0.2cm]
P\left[\sup_{t\le T}  | U\n_t | \ge \epsilon\right] &\le& \frac{4}{\epsilon^2} 
\EE \left[\int_0^t (Z\n_s)^2R_s\n ds\right] \le \frac{4C}{N\epsilon^2}, \label{eq:Doob}
\end{eqnarray}
where $C$ is a positive constant depending only on $\psi$. Thus $U\n_t\to 0$ almost surely as $N\to\infty$. This last property together with assumption (\ref{eq:init}) yield 
(\ref{C2}) and the announced (weak) relative compactness of the sequence $Z\n_t$. Hence, the sequence of probability measures $Q\n$, defined on  $\Dc_\mathcal{M}[0,T]$ and corresponding to the process $\mu\n_t$, is also relatively compact: this is a consequence of  classical projection theorems (see for instance Theorem 16.27 in \cite{Ka}). We are now in a position to state a  further important property.

Let  $Q$ the  limit point of some  arbitrary subsequence $Q^{(n_k)}$, as $n_k\to\infty$, and  $Z_t\egaldef \lim_{{n_k}\to\infty}Z^{(n_k)}_t$. Then  the support of $Q$ is a set of sample paths absolutely continuous with respect to the Lebesgue measure. Indeed, the application $\mu_t\to \sup_{t\le T}\log Z_t$ is continuous and we have the immediate bound
\[
\sup_{t\le T}\log Z_t \le \int_0^1 [ |\phi_a(x) | +  |\phi_b(x)|]dx,
\]
which holds for all $\psi_{a}, \psi_{b} \in \mathbf{C}^2[0,1]$. Hence, by weak convergence, any limit point $Z_t$ has the form
\begin{equation}\label{eq:lim}
Z_t[\phi_a,\phi_b] = \exp\Bigl[\int_0^1 [\rho(x,t)\phi_a(x)+(1-\rho(x,t)\phi_b(x) ] dx\Bigr],
\end{equation}
where  $\rho(x,t)$  denotes the limit density (a priori random) of the sequence of empirical measures  $\mu_t^{(m_k)}$ introduced in  (\ref{eq:emp1}).

\paragraph{P2 [A functional integral operator to characterize limit points]} 
This is somehow the Gordian knot of the problem. Relying on the above weak compactness property, our next result shows that any arbitrary limit point $Q$ is concentrated on a set of trajectories which are weak solutions of a \emph{functional integral equation} (IFE). 

First, by  (\ref{eq:martin1}), (\ref{eq:gen1}) and (\ref{eq:gen2}), we obtain at once 
\begin{equation} \label{eq:deriv1}
\begin{split}
&\frac{\partial (Z\n_t - U\n_t)}{\partial t}  = \\
& N^2\sum_{i\in\Gb\n}\widetilde{\lambda}_{ab}(i,N) \frac{\partial^2 Z\n_t}{\partial\phi_a(\frac{i}{N})\partial\phi_b(\frac{i+1}{N})} + \widetilde{\lambda}_{ba}(i,N) 
\frac{\partial^2 Z\n_t}{\partial\phi_a(\frac{i+1}{N}) \partial\phi_b(\frac{i}{N})}.
\end{split}
\end{equation}
It is worth remarking that  (\ref{eq:deriv1}) should be written, strictly speaking, as a  stochastic differential equation, which is well-defined since indeed all the underlying probability spaces emanate from a families of interacting Poisson processes.  

Replacing for a while  the quantities  $\phi_a(\frac{i}{N})$ and 
$\phi_b(\frac{i}{N})$ by variables $x\n_i$ and $y\n_i$ respectively, (\ref{eq:deriv1}) becomes
\begin{equation} \label{eq:deriv2}
 \frac{\partial (Z\n_t-U\n_t)}{\partial t} = N^2\sum_{i\in\Gb\n}
\alpha_{xy}(i,N) \frac{\partial^2 Z\n_t}{\partial x\n_i \partial y\n_{i+1}}
+ \alpha_{yx}(i,N)\frac{\partial^2 Z\n_t}{\partial y\n_i \partial x\n_{i+1}},
\end{equation}
where we have put
\begin{eqnarray*}
\alpha_{xy}(i,N) &= &\lambda_{ab}(N)\left[
 \exp\Bigl(\frac{ x\n_{i+1}  - x\n_i  + y\n_i  - y\n_{i+1}}{N} \Bigr) - 1\right], \\
 \alpha_{yx}(i,N) & = & \lambda_{ba}(N)\left[
 \exp\Bigl(\frac{ y\n_{i+1}  -y\n_i  + x\n_i  - x\n_{i+1}}{N} \Bigr) - 1\right].
 \end{eqnarray*}

We shall rewrite (\ref{eq:deriv2}) in the operator form
\begin{equation} \label{eq:deriv3}
 - \frac{\partial U\n_t}{\partial t}   = \Lc\n_t [Z\n_t],
\end{equation}
remarking  in the present setting that, for each finite $N$, $\Lc\n$ acts on the function space $\Cc^p{\bigl[-|\phi|,|\phi|]^{2\nn}}$, where
\begin{equation}\label{eq:norm}
|\phi| \egaldef \sup_{z\in[0,1]} \bigl(|\phi_a(z)|, |\phi_b(z)|\bigr),
\end{equation}
and $p$ is an  arbitrary positive number,  as $Z\n_t$ is analytic with respect to 
$\{\phi_a(.),\phi_b(.)\}$. The operator  $\Lc\n_t$ is  of parabolic type, but then in the wide sense, since here one can check  the quadratic form usually associated with the second order derivative terms is non definite, see e.g. \cite{EgSh}).

The key point will be to show that any limit point $\DD Z_t\stackrel{Law}{=}\lim_{n_{k}\to\infty}Z_t^{(n_k)}$ satisfies an IFE, obtained by studying the second order linear partial differential operators $\Lc\n_t$ along the sequence $n_k\to\infty$. 

To carry out the analysis of the limit sum coming in (\ref{eq:deriv2}) (which is a priori intricate), we propose a general approch, which aims at proving first that $Z_t$ is a \emph{weak solution} (or \emph{distributional} in the sense of Schwartz) of a Cauchy type operator. The line of argument will be sketched below. 

Beforehand, for the sake of shortness, it will be convenient to define the following cylinder sets, for $p=1,2\ldots$,
\[
\Uc^p_t \egaldef [-|\phi|,|\phi|]^p\times[0,t], \quad 
\Uc^p \egaldef [-|\phi|,|\phi|]^p .
\]

Introduce the operator $\widetilde{\Lc}\n_t$, which is the adjoint of  $\Lc\n_t$ in the Lagrange sense,  so that, for every function  $h\in\Cc^\infty_{0}(\Uc^{2\nn}_t)$,
\begin{eqnarray}\label{eq:opadj1}
\widetilde{\Lc}\n_t [h] &\egaldef&  
\frac{\partial h}{\partial t} + N^2\sum_{i\in\Gb\n}
 \frac{\partial^2 \bigl[\alpha_{xy}(i,N)h\bigr]}{\partial x\n_i \partial y\n_{i+1}}
+\frac{\partial^2 \bigl[\alpha_{yx}(i,N)h\bigr]}{\partial y\n_i \partial x\n_{i+1}} \nonumber \\
&\egaldef&  \frac{\partial h}{\partial t} + \Bc\n[h].
\end{eqnarray}

\begin{defin} [see e.g. \cite{Smi}, Part 1, section III.5]
A function $g\in L_2$ is said to be a weak (or distributional) solution of the Cauchy problem $\Lc\n g =0$ if, for all $h\in\Cc^\infty_{0}(\Uc^{2\nn}_T)$,
\begin{equation}\label{eq:opadj2}
\int_{\Uc^{2\nn}_T} g \widetilde{\Lc}\n_t[h] \,d\vec{u} \,dt = 0 ,
\end{equation}
where in the integral $\vec{u}$ denotes an arbitrary point in $\Uc^{2\nn}$.
\end{defin}

Multiplying equation (\ref{eq:deriv2}) by an arbitrary function  
$h\in\Cc^\infty_{0}(\Uc^{2\nn}_T)$, for fixed $T$ arbitrary positive, and then integrating twice by parts, we obtain, in agreement with (\ref{eq:opadj2}), 
\begin{equation}\label{eq:weak1}
\begin{split}
\int_{\Uc^{2\nn}_T} \left[  (Z\n_t - U\n_t)\frac{\partial h}{\partial t}  + 
 Z\n_t\Bc\n[h]\right]\, d\vec{u}\,dt
& = \\[0.2cm]
\int_{\Uc^{2\nn}} \bigl[(Z\n_T - U\n_T)h(\vec{u},T)-  Z\n_0h(\vec{u},0)\bigr]d\vec{u}\, .
\end{split}
\end{equation}

A brute force analysis of the adjoint operator could lead to a dead-end. A preliminary step will be to exploit carefully the estimates obtained in Lemma \ref{lem:esti}. This is the content of the next lemma.

\begin{lem}\label{lem:esti2}
The following partial differential equation holds.
\begin{equation} \label{eq:deriv4}
\begin{split}
\frac{\partial (Z\n_t-U\n_t)}{\partial t} & = 
\sum_{i\in\Gb\n} \mu \psi_{ab}' \Bigl(\frac{i}{N}\Bigr) \left[\frac{1}{2}\biggl( \frac{\partial Z\n_t}{\partial x\n_i} +  \frac{\partial Z\n_t}{\partial x\n_{i+1}}\biggr)-N \frac{\partial^2 Z\n_t}{\partial x\n_i \partial x\n_{i+1}} \right] \\[0.2cm]
& +\lambda \sum_{i\in\Gb\n} \psi_{ab}'' \Bigl(\frac{i}{N}\Bigr) 
\frac{\partial Z\n_t}{\partial x\n_{i+1}} + \Oc\Bigl(\frac{1}{N}\Bigr), 
 \end{split}
 \end{equation}
 where the term $\Oc\bigl(\frac{1}{N}\bigr)$ is in modulus uniformly bounded by 
 $\frac{C}{N}$, $C$ being a constant depending only on $\psi, \psi'$ and $\psi''$.  
 \end{lem} 
\begin{proof}
Immediate from equations (\ref{eq:estim1}) and (\ref{eq:estim2}).
\end{proof}

Starting from Lemma \ref{lem:esti2}, we will  present the two global guidelines of a functional approach, called  \textbf{G1} and \textbf{G2}. Basically, it relies on partial differential equations, whose variables are \emph{functions taken at points of the torus}. We think this might well extend to larger dimensions, although this assertion could certainly be debated.

\begin{itemize}
\item[ \textbf{G1}] Intermediate reduction to an almost sure convergence context. This can be achieved by means of the extended Skohorod coupling theorem (see Corollary 6.12  in \cite{Ka}), which in brief says that, if a sequence of real random variables $(\xi_k)$ is such that $\lim_{k\to\infty}f_k(\xi_k)= f(\xi)$  converges in distribution, then there exist a probability space  $\mathcal{V}$ and  a new random sequence $\widetilde{\xi_k}$, such that  $\widetilde{\xi_k}\stackrel{\Lc}{=}\xi_k$ and $\lim_{k\to\infty} f_k(\widetilde{\xi_k})= f(\xi)$, almost surely in $\mathcal{V}$, with $\widetilde{\xi}\stackrel{\Lc}{=}\xi$. Here this theorem will be applied to the family   $Z^{(n_{k})}_t$,  which thus gives rise a new sequence  denoted by $Y^{(n_{k})}_t$  in the sequel. This step is in no way obligatory, but just a matter of taste. Indeed, one could still keep on with weak convergence context and use Alexandrov's portmanteau theorem (see e.g. \cite{EtKu})
 whenever needed.
\item[ \textbf{G2} ] For each finite $N$, we can consider the quantities 
\[
 \psi_{ab}' \bigl(\frac{i}{N}\bigr), \psi_{ab}'' \bigl(\frac{i}{N}\bigr), \ i = 1,\dots,N,
\]
as \emph{constant parameters, while the  $x\n_i$'s are free variables}. This  is clearly feasible, choosing for instance $\phi_a(.), \phi_b(.)$ in the class of polynomials of degree at least  $3N$. Also, from now on, the functions 
$\phi_a$ and  $\phi_b$ will be supposed to belong to  $\Cc^\infty_{0}(K)$, for some compact $K\in\Rc$ containing the interval $[0,1]$.
\end{itemize}

Then, according to  \textbf{G1}, we rewrite  (\ref{eq:deriv4}) as
\begin{equation}\label{eq:esti3}
-\frac{\partial U\n_t}{\partial t} \egaldef \Ac\n_t[Y\n_t] + \Oc\Bigl(\frac{1}{N}\Bigr),
\end{equation}
where $\Ac\n_t$ is viewed as an operator of parabolic type with \emph{constant} coefficients and  domain    $\Cc^\infty_{0}(\Uc^{\n}_T)$.
\[
\begin{split}
\Ac\n_t[g] & \egaldef -\frac{\partial g}{\partial t} +\sum_{i\in\Gb\n} \mu \psi_{ab}' \Bigl(\frac{i}{N}\Bigr) \left[\frac{1}{2}\biggl( \frac{\partial g}{\partial x\n_i} + \frac{\partial g}{\partial x\n_{i+1}}\biggr) -N \frac{\partial^2 g}{\partial x\n_i \partial x\n_{i+1}} \right] \\[0.2cm]
& +\lambda \sum_{i\in\Gb\n} \psi_{ab}'' \Bigl(\frac{i}{N}\Bigr) 
\frac{\partial g}{\partial x\n_{i+1}},
\end{split}
\]
remembering that $\psi_{ab}=\phi_{a} - \phi_{b}$. The term $\Oc\Bigl(\frac{1}{N}\Bigr)$
 in (\ref{eq:esti3}) stands for an operator having  a negligible range for $N\to\infty$.

Let $\widetilde{\Ac}\n_t[h]$  denote the adjoint of  $\Ac\n_t$. Then
\begin{equation}\label{eq:opadj3}
\begin{split}
\widetilde{\Ac}\n_t[h]  & = \frac{\partial h}{\partial t} - \sum_{i\in\Gb\n} \mu \psi_{ab}' 
\Bigl(\frac{i}{N}\Bigr) \left[ \frac{1}{2}\biggl( \frac{\partial h}{\partial x\n_i} + 
\frac{\partial h}{\partial x\n_{i+1}}\biggr) + N \frac{\partial^2 h}{\partial x\n_i \partial x\n_{i+1}}\right]  \\[0.2cm] 
& -\lambda \sum_{i\in\Gb\n} \psi_{ab}'' \Bigl(\frac{i}{N}\Bigr) 
\frac{\partial h}{\partial x\n_{i+1}},
\end{split}
\end{equation}
and, for any $h\n\in\Cc^\infty_{0}(\Uc\n_T)$,  we have
\begin{equation} \label{eq:opadj4}
\begin{split}
 \int_{\Uc\n} \left[(Y\n_T - U\n_T)h\n(\vec{u},T) - Y\n_0h\n(\vec{u},0)\right]
  \delta\vec{u} & =  \\[0.2cm]
 \int_{\Uc\n_t}  (Y\n_t - U\n_t)\widetilde{\Ac}\n_t[h\n]\,\delta\vec{u}\,dt  + 
\Oc\Bigl(\frac{1}{N}\Bigr) \,,
  \end{split} 
  \end{equation}
where $\delta\vec{u}$ in (\ref{eq:opadj4}) represents the differential volume element 
\[
\DD \delta\vec{u} = dx\n_1dx\n_2...dx\n_N= \delta\phi_a\Bigl(\frac{1}{N}\Bigr) \delta\phi_a\Bigl(\frac{2}{N}\Bigr)\ldots \delta\phi_a(1).
 \]

The next step is to make a suitable choice of the function $h$ in  (\ref{eq:opadj4}) in order to extract a meaningful information on the limit operator, as $N\to\infty$.
Keeping in mind that  the random variables  $Y\n_t$ are defined on the implicit probability space $\mathcal{V}$ introduced in \textbf{G1}, we state the following result.

 \begin{lem} \label{lem:esti3}  \mbox{\emph{\textbf{[Claim1]}}}
Let $\DD Y_t = \lim_{n_k\to\infty} Y^{(n_k)}_t \ a.s.$,  and let $W\in\Cc^\infty_{0}(K)$ be a fixed compact space. Then,  for a class of test functions $k$ properly chosen in a subset of $\Cc^\infty_0(W\times[0,T])$, we have
\begin{equation} \label{eq:opfunc}
\begin{split} 
& \int_0^T dt \int Y_t [\phi(.)]\, F_t [\phi(.)]  \delta\phi(.) = \\[0.2cm]
& \int \bigl[k\bigl(\phi(.),T\bigr) Y_T[\phi(.)] - k\bigl(\phi(.),0\bigr) Y_0[\phi(.)] \bigr] \delta\phi(.)\, ,
\end{split}
\end{equation}
where 
\[
\begin{split}
 F_t [\phi(.)] & =  \frac{\partial k\bigl(\phi(.),t\bigr)}{\partial t} \\[0.2cm]
& -\int_0^1\left[\mu\frac{\partial^2k\bigl(\phi(.),t\bigr)}{\partial\phi^2(x)} \psi'(x) + 
\bigl[\mu\psi'(x)+\lambda\psi''(x)\bigr]\frac{\partial k\bigl(\phi(.),t\bigr)}{\partial\phi(x)} \right] dx .
\end{split}
\]
Moreover (\ref{eq:opfunc}) yields the Cauchy problem posed in  (\ref{eq:Cauchy}).
\end{lem}

 \begin{proof}  
 We only sketch the main  lines of argument.
 \begin{itemize}
 \item[$\bullet$] The almost sure convergence of $Y_t\n$  and  $U_t\n$, respectively to $Y_t$ and $0$ [by (\ref{eq:Doob})], will yield (\ref{eq:opfunc}) provided that in  (\ref{eq:opadj4}) the functions $h\n$ are properly chosen, to ensure the existence of the limit sums coming in (\ref{eq:opadj3}), as $N\to\infty$. 
 
 Setting  $\vec{x}\n \egaldef (x\n_1,x\n_2,\ldots,x\n_N)$, 
 he reader can convince himself that it suffices to take in (\ref{eq:opadj4})
 \[
 h\n\Bigl(\frac{\vec{x}\n}{N}, t\Bigr) = k(\vec{x}\n,t),
 \]
 where $k\in\Cc^\infty_0(W\times[0,T])$,  in which case convergent Riemann sums are obtained in (\ref{eq:opadj3}).
 \item[$\bullet$] From  Skohorod's coupling theorem,  $Y_t$ does satisfy an equation of the form (\ref{eq:lim}). Hence, we can write  the following functional derivatives (which are plainly of a Radon-Nykodym nature)
\[
\begin{cases}
\DD \frac{\partial Y_t}{\partial\phi(.)} = \rho(.,t) Y_t  ,\\[0.3cm]
\DD \frac{\partial^2 Y_t}{\partial\phi^2(.)} = \rho^2(.,t)Y_t .
\end{cases}
 \]
 \item[$\bullet$] To derive  (\ref{eq:Cauchy}),  one has to pick out $k$ from a class of  convolution test functions, depending  on some parameter $\epsilon$ and properly converging in the space of Schwartz distributions.
\end{itemize}
\end{proof} 
\paragraph{P3 [Uniqueness]} \mbox{ }
The problem of uniqueness  of weak solutions of the Cauchy problem (\ref{eq:Cauchy}) for nonlinear parabolic quation is in fact already solved in  the literature. We refer the reader for instance to \cite{EgSh} for a wide bibliography on the subject. 

To conclude the analysis of Theorem \ref{theo:main}, it suffices to switch back to 
$Z\n _t[\phi_a,\phi_b]$, which converges in distribution to $Z_t[\phi_a,\phi_b]$. Hence, the random measure $\mu\n_{t}$ converges in distribution to a deterministic measure, which is a peculiar situation impling also convergence in probability.
\end{proof}

\section{The n-species model}

We  will state a conjecture about  hydrodynamic equations for the $n$-species model, briefly introduced in section  \ref{sec:clock}, in the so-called \emph{equidiffusion} case,
 precisely defined hereafter.
 
 The $n$-species system  is said to be \emph{equidiffusive} whenever there exists a constant 
 $D$, such that, for all pairs $(k,l)$,
 \[
 \lim_{N\to\infty} \frac{\lambda_{kl}(N)}{N^2} = D.
 \]
Letting
\[
\alpha_{kl} \egaldef  \lim_{N\to\infty} \log\frac{\lambda_{kl}(N)}{\lambda_{lk}(N)}, 
\]
we assert that the following hydrodynamic system holds \textbf{[Claim2]}.
\[
\frac{\partial\rho_k}{\partial t} = D\left[\frac{\partial^2\rho_k}{\partial x^2}
+\frac{\partial}{\partial x}\Bigl(\sum_{l\ne k}\alpha^{lk}\rho_k\rho_l\Bigr)\right], \ k=1,\ldots,n.
\]
The idea is  to apply the functional approach  presented in this paper. 

\nocite*
\bibliography{refer1}

\begin{thebibliography}{10}

\bibitem{ArHeRi}
{\sc P.~Arndt, T.~Heinzel, and V.~Rittenberg}, {\em Stochastic models on a ring
  and quadratic algebras. the three-species diffusion problem}, J. Phys. A:
  Math. Gen., 31 (1998), pp.~833--843.

\bibitem{Berge}
{\sc C.~Berge}, {\em Th\'eorie des Graphes et ses Applications}, vol.~II of
  Collection Universitaire des Mathématiques, Dunod, 2~ed., 1967.

\bibitem{BeLa}
{\sc L.~Bertini, A.~De~Sole, D.~Gabrielli, G.~Jona~Lasinio, and C.~Landim},
  {\em Current fluctuations in stochastic lattice gases}, Phys. Rev. Lett., 94
  (2005), p.~030601.

\bibitem{Bil}
{\sc P.~Billingsley}, {\em Convergence of Probability Measures}, Wiley Series
  in Probability and Statistics, John Wiley \& Sons Inc., 2~ed., 1999.

\bibitem{Bu}
{\sc J.~M. Burgers}, {\em A mathematical model illustrating the theory of
  turbulences}, Adv. Appl. Mech., 1 (1948), pp.~171--199.

\bibitem{ClDeEv}
{\sc M.~Clincy, B.~Derrida, and M.~Evans}, {\em Phase transition in the {ABC}
  model}, Phys. Rev. E, 67 (2003), pp.~6115--6133.

\bibitem{MaPr}
{\sc A.~De~Masi and E.~Presutti}, {\em Mathematical Methods for Hydrodynamic
  Limits}, vol.~1501 of Lecture Notes in Mathematics, Springer-Verlag, 1991.

\bibitem{DeEvHaPa}
{\sc B.~Derrida, M.~Evans, V.~Hakim, and V.~Pasquier}, {\em Exact solution for
  1d asymmetric exclusion model using a matrix formulation}, J. Phys. A: Math.
  Gen., 26 (1993), pp.~1493--1517.

\bibitem{EgSh}
{\sc Y.~Egorov and M.~Shubin}, eds., {\em Partial Differential Equations},
  vol.~I-II-III of Encyclopedia of Mathematical Sciences, Springer Verlag,
  1992.

\bibitem{EtKu}
{\sc S.~Ethier and T.~Kurtz}, {\em Markov Processes, Characterization and
  Convergence}, John Wiley {\&} Sons, 1986.

\bibitem{EvFoGoMu}
{\sc M.~Evans, D.~P. Foster, C.~Godrèche, and D.~Mukamel}, {\em Spontaneous
  symmetry breaking in a one dimensional driven diffusive system}, Phys. Rev.
  Lett., 74 (1995), pp.~208--211.

\bibitem{EvKaKoMu}
{\sc M.~Evans, Y.~Kafri, M.~Koduvely, and D.~Mukamel}, {\em Phase {S}eparation
  and {C}oarsening in one-{D}imensional {D}riven {D}iffusive {S}ystems}, Phys.
  Rev. E., 58 (1998), p.~2764.

\bibitem{FaFu}
{\sc G.~Fayolle and C.~Furtlehner}, {\em Dynamical {W}indings of {R}andom
  {W}alks and {E}xclusion {M}odels. {P}art {I}: Thermodynamic limit in $\mathbb
  {Z}^2$}, Journal of Statistical Physics, 114 (2004), pp.~229--260.

\bibitem{FaFu2}
{\sc G.~Fayolle and C.~Furtlehner}, {\em Stochastic deformations of sample
  paths of random walks and exclusion models}, in Mathematics and computer
  science. III, Trends Math., Birkh\"auser, Basel, 2004, pp.~415--428.

\bibitem{GoLu}
{\sc C.~Godr\`eche and J.~Luck}, {\em Nonequilibrium dynamics of urns models},
  J. Phys. Cond. Matter, 14 (2002), p.~1601.

\bibitem{Ka}
{\sc O.~Kallenberg}, {\em Foundations of Modern Probability}, Springer, 2~ed.,
  2001.

\bibitem{KPZ}
{\sc M.~Kardar, G.~Parisi, and Y.~Zhang}, {\em Dynamic scaling of growing
  interfaces}, Phys. Rev. Lett., 56 (1986), pp.~889--892.

\bibitem{Kel}
{\sc F.~P. Kelly}, {\em Reversibility and stochastic networks}, John Wiley \&
  Sons Ltd., 1979.
\newblock Wiley Series in Probability and Mathematical Statistics.

\bibitem{KiLa}
{\sc C.~Kipnis and C.~Landim}, {\em Scaling limits of Interacting Particles
  Systems}, Springer-Verlag, 1999.

\bibitem{LaBaRa}
{\sc R.~Lahiri, M.~Barma, and S.~Ramaswamy}, {\em Strong phase separation in a
  model of sedimenting lattices}, Phys. Rev. E, 61 (2000), pp.~1648--1658.

\bibitem{Li}
{\sc T.~M. Liggett}, {\em Stochastic Interacting Systems: Contact, Voter and
  Exclusion Processes}, vol.~324 of Grundlehren der mathematischen
  {W}issenschaften, Springer, 1999.

\bibitem{MUR}
{\sc J.~Murray}, {\em Mathematical Biology}, vol.~19 of Biomathematics,
  Springer-Verlag, 2~ed., 1993.

\bibitem{RUD}
{\sc W.~Rudin}, {\em Functional Analysis}, International Series in Pure and
  Applied Mathematics, McGraw-Hill, 2~ed., 1991.

\bibitem{Smi}
{\sc V.~Smirnov}, {\em Cours de Mathématiques Supérieures}, vol.~IV, \'Editions
  MIR - Moscou, traduction fran\c caise, 1984.

\bibitem{Sp}
{\sc H.~Spohn}, {\em Large Scale Dynamics of Interacting Particles}, Springer,
  1991.

\end{thebibliography}
\bibliographystyle{siam}

\end{document}